\documentclass[conference]{IEEEtran}
\IEEEoverridecommandlockouts
\usepackage{cite}
\usepackage{amsmath,amssymb,amsfonts}
\usepackage{algorithmic}
\usepackage{graphicx}
\usepackage{textcomp}
\usepackage{xcolor}
\usepackage{url} 
\def\BibTeX{{\rm B\kern-.05em{\sc i\kern-.025em b}\kern-.08em
    T\kern-.1667em\lower.7ex\hbox{E}\kern-.125emX}}
\usepackage{float}
\usepackage[font={footnotesize}]{caption}

\makeatletter
\def\ps@IEEEtitlepagestyle{%
  \def\@oddfoot{\mycopyrightnotice}%
  \def\@evenfoot{}%
}
\def\mycopyrightnotice{
  {\footnotesize 979-8-3503-8639-4/24/\$31.00 ©2024 IEEE\hfill} 
  \gdef\mycopyrightnotice{}% just in case
}
\begin{document}

\title{Smart Health Software to Support Rescue Personnel in Emergency Situations}

\author{\IEEEauthorblockN{Abu Shad Ahammed}
\IEEEauthorblockA{\textit{Chair of Embedded Systems} \\
\textit{University of Siegen}\\
Siegen, Germany \\
ORCID: 0009-0007-6715-8098}
\and
\IEEEauthorblockN{Roman Obermaisser}
\IEEEauthorblockA{\textit{Chair of Embedded Systems} \\
\textit{University of Siegen}\\
Siegen, Germany \\
ORCID: 0009-0002-4483-1503}
}

\maketitle

\begin{abstract}
Rescue stations around the world receive millions of emergency rescue calls each year, most of which are due to health complications. Due to the high frequency and necessity of rescue services, there is always an increasing demand for quick, accurate, and coordinated responses from rescue personnel to save lives and mitigate damage. This paper introduces a rescue health management software solution designed to improve the efficiency and effectiveness of rescue situational awareness by rapidly assessing the health status of emergency patients using AI-driven decision support systems. The novelty in this software approach is it's user-centered design principles to ensure that its solutions are specifically tailored to meet the unique requirements of emergency responders. It used pre-trained machine learning models with rescue data and accepted new patient's input data to provide a probability of the major health complications so that rescue personnel can expedite treatment plan following the outcome. The paper focuses primarily on the software development and implementation steps with three use cases, while also providing a short overview of the previous machine learning-based development phases.
\end{abstract}


\section{Introduction}
Rescue emergencies are generally quite strenuous and challenging because they deal with human lives in a situation where it is difficult to accurately assess the health status of a distressed patient due to personal and work limitations. Often rescue personnel face high levels of mental and physical stress, trauma, and emotional stress during rescue operations, especially when dealing with injured and vulnerable victims. This certainly impacts their decision-making abilities and overall well-being.  Time, the main constraint in such rescue situations, also plays a vital role in decision making.  Therefore, it is necessary to recognize relevant situations, i.e., health complications of the rescue patients on site, and to take appropriate first aid measures. To help rescue personnel identify the health complication of a patient in such stressful situations, a machine learning (ML) based software solution can be provided that can save time and improve care management. The principal advantage of using ML in healthcare is its ability to analyze large amounts of data in a short period of time to provide a quick diagnosis.\\
Artificial intelligence techniques ranging from machine learning to deep learning are widely used in health care for disease diagnosis, treatment path recommendation, and patient health management. Machine learning (ML) algorithms, for example, support vector machine (SVM), K-nearest neighbors (KNN), extreme gradient boosting (XGB), etc. are widely used in the prognosis of diseases like skin cancer, diabetes, epilepsy, and others\cite{kumar2023artificial}. The purpose of using artificial intelligence in the modern health sector is to reduce the workload of physicians, minimize errors and diagnostic times, and improve performance in the prediction and detection of various diseases. In this paper, we discuss a software solution developed with multiple ML algorithms that can provide a quality health diagnosis so that the treatment path can be expedited. The data used to develop the AI algorithms were obtained from the Rescue Station in Siegen-Wittgenstein. Although the detection precision achieved through AI algorithms is quite significant, it should be mentioned that the experience of physicians in diagnosing such complications is undeniable\cite{ahammed2023time}. The interpretability of AI-driven disease recognition systems can pose a challenge, as the algorithms may lack transparency in decision-making processes, making it difficult for healthcare providers to understand and trust algorithmic output\cite{mooghali2023barriers}. Therefore, the software developed in our study is intended to assist rescue teams in their decision-making process, rather than serving as the primary solution. The focus was to design reliable and robust software that can consider dynamic changes in patient health data and support decision-making.\\
The structure of the rest of the manuscript is as follows. Section \ref{lit} presents a meta-analysis of the current research, methods, and techniques by various authors. Section \ref{data} covers the organization and management of the raw rescue data obtained from the rescue station. Section \ref{ml} briefly discusses the machine learning algorithms employed to create the detection model. The significance of health management software is discussed in Section \ref{sw}. The section also explains the requirement analysis of the software developed in our research. The design steps of the software and its implementation process are analyzed in Section \ref{impl}. To assess the impact and performance of the software, Section \ref{test} provides examples of use cases. Finally, Section \ref{conc} concludes with a discussion of the existing challenges and future possibilities.
\section{Literature Review}\label{lit}
Health management software to support rescue personnel identify the correct health complications is unique research and novel in this field. From our exploration, we have found plenty of medical research on developing ML-based software for disease identification and care management. However, we discovered that although machine learning has been successfully applied in various scientific domains of health management, including the detection of health complications, its importance in aiding emergency rescue patients remains under explored. The type of data that we are utilizing to develop the software, namely rescue data, is unprecedented and novel work in this research field.\\
We found some important studies highlighting the potential of machine learning algorithms in predicting complications, mortality rates, length of hospital stays, and overall health improvement \cite{shillan2019use, bertini2022using, moraru2020artificial}. The researchers have used ML algorithms in different healthcare specialties, including ophthalmology, dentistry, and prosthodontics, demonstrating their versatility and effectiveness in care management. Machine learning has been found to be particularly significant in the diagnosis of metabolic diseases like diabetes or in the prediction of stroke based on early patient analysis and the integration of various data sources like clinical symptoms, laboratory tests, genetic information, and imaging results \cite{mei2020artificial, brasil2019artificial}. Our previous studies \cite{ahammed2022detection, ahammed2022novel} on predicting cardiovascular and respiratory complications in emergency patients demonstrated a precision that exceeded 90\% with XGB as it was found to be the model with the best performance. There are other significant healthcare studies similar to ours showing ML, in particular, the XGB model demonstrating an accuracy greater than 90\% for neurological or cardiovascular diseases \cite{nabi2021efficient,nagavelli2022machine}.\\
The significance of health management software is well recognized with the advancement of deep learning algorithms. Many of these software are capable of analyzing electroencephalography or electrocardiography signals, radiography, X-rays, or CT scan images to detect the actual condition of a patient \cite{shoeb2010application, abbasi2023automatic}. Yang et al. \cite{yang2023smart} had provided a conceptual overview on how integration of machine learning and deep learning-based health care system promises to deliver precise prognostic predictions for diseases like COVID-19, thus improving physician decision making.
\section{Rescue Data Management}\label{data}
The data used in this research were collected through a KMU-Innovativ research project-`KIRETT' from Germany where the aim was to detect rescue situations to expedite treatment path of patients. The Siegen-Wittgenstein rescue station provided these data for research purposes, which was recorded from January 2012 to August 2021 in various rescue scenarios in the city of Siegen, Germany. The data set contains anonymous information on 273,183 previous rescue events, each identified with various health complications and comprising more than 452 unique attributes. The rescue cases represent the history of rescue operations for different patients suffering from various complications. Each case was recorded with information-like rescue protocols used, the location of the rescue events, the clinical record of the patients, the health prognosis, vital signs, the initial impression of the rescuers and the treatment performed.\\
The data from the rescue station was recorded in an emergency environment where the collection of information is quite difficult and sometimes gets distorted. Hence, it was not surprising to find misspelled text information, and abnormal values of health vitals like blood glucose, respiratory rate, pulse rate, etc. It initiated the need to perform data preprocessing techniques to make the data usable for machine learning applications. The following preprocessing steps were considered to make the raw data usable:
\begin{itemize}
    \item \textbf{Data Integration}: More than 80 rescue files were obtained from the rescue station, each containing unique information about patients involved in similar rescue incidents. Using the case number as the data key, these files were merged to make a unique database containing all the information from a single rescue event.
    \item \textbf{Data Reduction}: The rescue data set contained 452 columns, in which many of those columns did not contain any significant information. Some columns had duplicate statistics, which were removed or unified. Finally, 380 columns were left to be used for the next process.
    \item \textbf{Data Filtering}: There was a lot of misspelled and ambiguous information in the database which was filtered through multiple steps. Firstly, the health vital information was checked through the Interquartile Range(IQR) using python. It identified all rescue cases with outliers that were later replaced with the average value of the respective health vital calculated from the database. Correcting textual information was challenging due to the size of the data and identifying whether the data should fall under a misspelled category or acronyms used due to a lack of time. Eventually, parsing the data was not the only solution. So, Microsoft Excel data processing techniques like filtering and sorting and natural language tool kits in Python were utilized to unify similar words. 
\end{itemize}
\section{Algorithm Development}\label{ml}
Machine learning algorithms are a state-of-the-art technology used to aid medical professionals in the early detection of diseases or identifying health complications to assess the risk and prognosis of the disease \cite{dai2015prediction, nadeem2021fusion}. The health management software presented in this paper is also developed with ML as it's core in diagnosis of the major health complications observed in rescue record: Cardiovascular, Respiratory, Neurological, Psychiatric, Abdominal, and Metabolic. As machine learning is not the main focus of this paper, only a short overview of the development steps are given below.
\begin{itemize}
    \item \textbf{Feature Selection}: The plan was to develop individual ML models for each major health complication. So, for each ML model, we had to identify compatible features to prepare the training data. To select the initial features of a health complication detection model, we discussed with medical professionals and conducted our own research to identify which of these features are related to a health complication. Not all relevant features were available in the data set. In addition, some features were extracted from textual information using natural language processing tool kits and later converted to categorical features.
    \item \textbf{Selection of Algorithms}: Based on our findings on state-of-the-art research \cite{abas2024machine, liu2022machine}, we selected seven algorithms like Extreme gradient boosting (XGB), Support vector machine (SVM), Random forest (RF), K-nearest neighbor (KNN), Naive Bayes (NB), Logistic regression (LR), Artificial neural network (ANN) with 2 hidden layers. All models were developed as a binary prediction model.
    \item \textbf{Feature Optimization for the Individual Model}: Not all features are suitable for a machine learning algorithm in terms of achieving the highest accuracy. So, the recursive feature elimination with cross-validation (RFECV) technique was applied to identify these trivial features and remove them. For the algorithms where RFECV is ineffective, we checked the feature importance score and correlation values to eliminate unnecessary features. The finally selected features contained information as follows:
    \begin{itemize}
        \item Health vitals: Blood pressure, Respiratory rate, Circulation state, Mean arterial pressure, Body temperature, SpO2, Blood glucose, Pulse rate, Glasgow comma scale value
        \item Pre-illness information
        \item Health abnormality symptoms pertained to certain complication
        \item Neurological disturbances
        \item Patient's intoxication with drug or alcohol
    \end{itemize} 
    \item \textbf{Hyperparameters Tuning}: To find a set of optimal hyperparameters that can tweak the model performance, we chose `Grid Search' and `Random Search' tuning algorithms. For ANN algorithms, we also used 'Hyperband' and 'Bayesian Optimization' and then passed combination of specified hyperparameters to choose the best one. The purpose of tuning was to maximize the model's performance without overfitting, underfitting, or creating a high variance. 
    \item \textbf{Model Training and Evaluation}: The next step was to train the ML models for each complication to find the model that delivers the best performance. Individual training data for each complication were prepared with output labels 1 and 0. 1 indicates whether a feature vector in the training data corresponds to patients,  and 0 for non-patients. The evaluation criterion considered to judge their performance was precision, accuracy, and recall. After testing the algorithms with an unknown test set, we found that, as expected, the XGB model delivers the best performance for all six complications. In Table \ref{tab:perf}, the detailed results of the XGB models for each complication are provided.
\end{itemize}
\begin{table}[!t]
\renewcommand{\arraystretch}{1.3}
\caption{\textsc{Performance of XGB Model for Six Major Health Complications}}
\label{tab:perf}
\centering
\begin{tabular}{|l|l|l|l|}
\hline
\multicolumn{1}{|c|}{\textbf{\begin{tabular}[c]{@{}c@{}}Health\\ Complicatrion\end{tabular}}} & \textbf{Precision} & \textbf{Accuracy} & \textbf{Recall} \\ \hline
Cardiovascular                                                                                & 98.51\%            & 96.66\%           & 94.22\%         \\ \hline
Respiratory                                                                                   & 93\%               & 92.55\%           & 91.57\%         \\ \hline
Neurological                                                                                  & 84.2\%             & 85.93\%           & 87.99\%         \\ \hline
Psychiatric                                                                                   & 89.63\%            & 87.81\%           & 87.07\%         \\ \hline
Abdomen                                                                                       & 93.48\%            & 93.42\%           & 94.34\%         \\ \hline
Metabolic                                                                                     & 86.41\%            & 86.86\%           & 81.04\%         \\ \hline
\end{tabular}
\end{table}
Although XGB model displays the best performance, we also considered deep learning model like artificial neural network in our software due to its capability of modelling biological system with big data and fast learning of new data \cite{saharuddin2020constitutive}.
\section{Rescue Health Management Software}\label{sw}
Rescue health management software is a subset of health management software that focuses on supporting rescue personnel to identify a patient's health situation and expedite the treatment path. The health management software used by clinical institutes or hospitals is nothing more than a combination of specialized digital tools or platforms designed to improve various aspects of healthcare delivery, such as patient care, administrative tasks, and data management. Currently, many medical organizations are also using it as a clinical decision support system to support healthcare providers in quickly accessing patient information, monitoring treatment plans, and tracking health outcomes, which ultimately leads to improved patient care and treatment results \cite{ferrando2015effective}. The software developed in this research is expected to play an important role in emergency events by providing a comprehensive decision support system that can minimize the risk of errors associated with manual decision-making for situation recognition.
\subsection{Requirement Analysis}\label{req}
Software requirements are specifications that define the functionalities, constraints, and attributes of a software system. According to IEEE, the requirements are like feature or non-functional constraint that the software must provide to fulfill the needs of users and other stakeholders \cite{ieee_software_requirements}. Some of the common software requirements we found relevant to our software are:
\begin{itemize}
    \item \textbf{Functional Requirements}: Identify the components and methods that are necessary for the software. The requirement must describe the interactions between the system and its environment, independently of their implementation process. 
    \item \textbf{Non-Functional Requirements}: These requirements specify the characteristics of the system, including performance, usability, reliability, flexibility, and scalability. The performance of the software can be verified with two criterion: response time and the ability to interpret the outcome of machine learning algorithms apprehensively.
    \item \textbf{System Requirements}: These requirements are to ensure that the software is compatible with smart devices and operating systems used by rescue personnel.
    \item \textbf{User Requirements}: The most important requirement is intended for specific users, i.e., rescue personnel. The software should be designed considering rescue situations, the rescuers' work procedure, and the environment. It should not add additional responsibility; instead, it should shorten the work load with quality improvement.
\end{itemize}
\section{Software Design and Implementation}\label{impl}
\begin{figure}[!b]
%\captionsetup{justification=left}
\includegraphics[scale=0.18]{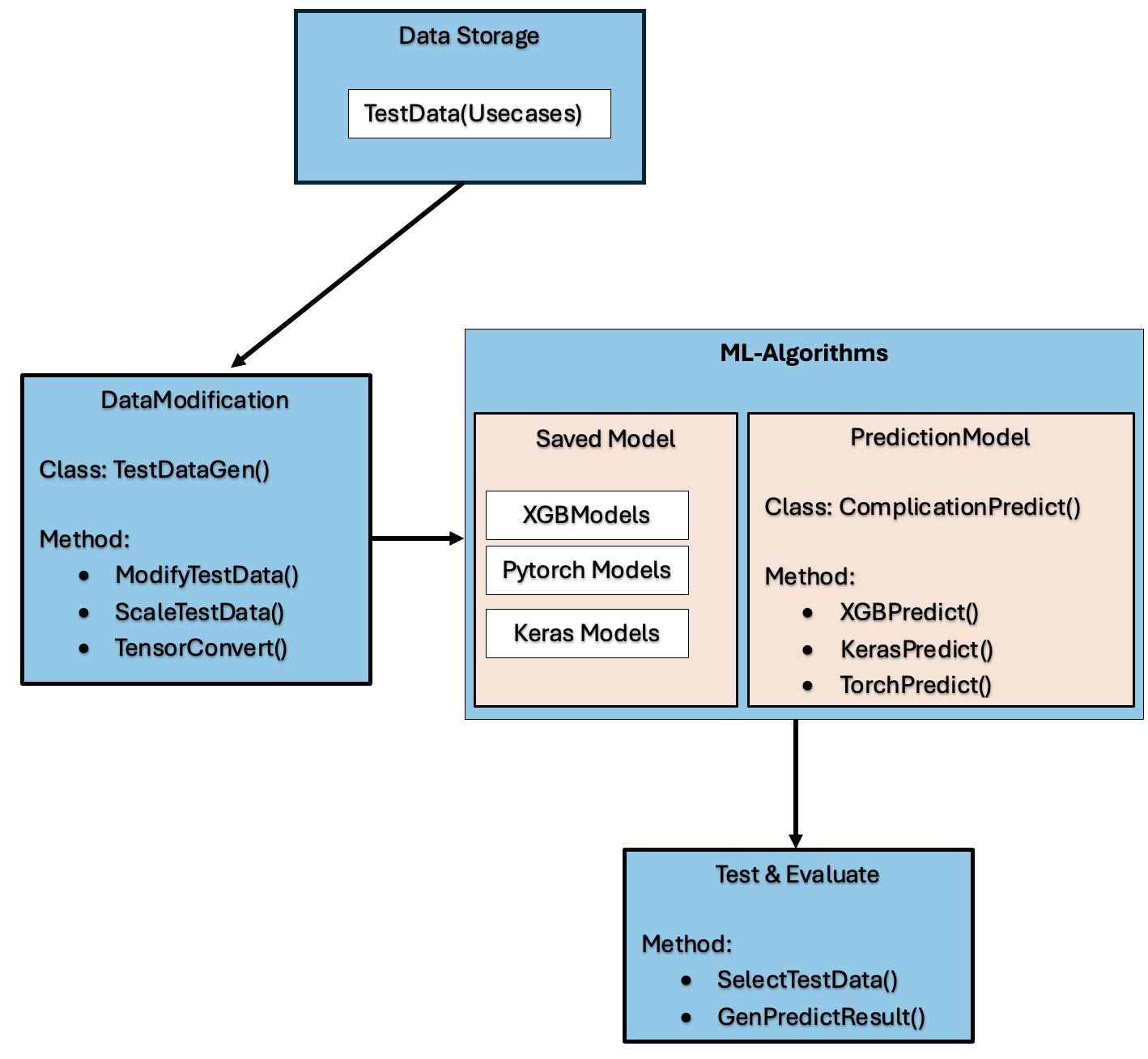}
\DeclareGraphicsExtensions{.png,.pdf}
\caption{Rescue Health Management SW Architecture to Predict Probability of Health Complications}
\label{Fig:arch}
\end{figure}

The software was developed using Python version 3.11 and Pycharm interactive development environment(IDE) in windows operating system. The Python object-oriented programming paradigm was utilized, which is basically a computer programming design philosophy or methodology that organizes or models software architecture around data or objects, rather than functions and logic \cite{pun1990design}. Python project template cookiecutter was used and customized to design the high-level template of the software architecture \cite{cookiecutter}. For low-level design, the classes and methods intended to integrate the ML algorithms and required test data were outlined in a virtual environment. The complete software architecture can be seen in Fig. \ref{Fig:arch}.
\\
The SW architecture is defined with four main blocks. The data storage block contains six test case data with each having 32 health characteristics of the patients. The features represent patient's health vitals, pre-illness information, injury or trauma, any communication or consciousness disorder, or presence of complication-specific abnormal signs. As the test data cannot be used directly to test the algorithms rather need to be preprocessed like converting to tensors for PyTorch model or scaling, the Data Modification block was designed. A method for modifying the health vital information of the test data is also defined in this block. The method will adjust the health vital according to user preferences and evaluate the impact of these changes on the performance of the ML models.
ML-algorithms block saves all the models i.e. ANN models and best performing XGB models that have been created for each complication. The ANN models were created using both the Keras and PyTorch frameworks to flexibly support future deployment at different hardware platforms. Another sub-block in that block was the 'PredictionModel' which is the soul of the software and runs the major operations i.e., predicting test data modified from the DataModification block using the savel ML algorithms. The Test \& Evaluate block contains a software test script and is intended to present the prediction results in a structured manner, making it easier to comprehend the software output and decide on subsequent steps. When test data are provided as input to the ML algorithms, this block displays the probability percentage of the six complications derived from the algorithms.
If a user wants to implement the software, the following procedures are required to be performed:
\begin{itemize}
    \item Run the file 'Demonstration.Py' through Python terminal. It will ask for the test file where the use cases were saved. 
    \item Select the test file and a test case
    \item An option will be provided to the user if a change in test data is intended. The user needs to type yes or no based on his requirement.
    \item If yes was selected in the previous step, the program will ask how much deviation in percentage is intended
    \item Finally, the software will show the percentage probability of each complication for the selected test data. The probability outcome will be displayed for both Keras and Pytorch models. If modification was done, then the results based on modification will also be displayed in parallel.
\end{itemize}
During the software design phase, after defining the classes and methods, each of them was carefully tested so that any flaw can be identified at the beginning. Detecting issues in the primary stage helps developers minimize the chances that these issues will be affected in later phases. Furthermore, tackling and fixing problems early is simpler, more economical, and reduces the need for redundant work \cite{boehm2001defect, westland2002cost}. 
\section{Software Testing and Evaluation}\label{test}
Software testing is a significant approach in the software development lifecycle that ensures the successful functionality of a software program. It involves verifying and validating that a software application is free of bugs and meets the functional, non-functional, and user requirements set at the planning phase without any contrast. For the proposed health management software, a testing environment is established to assess its performance according to the requirement analysis detailed in Section \ref{req}. In this study, three randomly chosen use cases will serve as test data to assess the effectiveness of the software. The data of the use cases, as can be seen in Table \ref{tab:ucd}, were extracted from the original rescue events with the diagnostic information of the rescue personnel and not used during the training of the ML models.
\begin{table}[H]
\renewcommand{\arraystretch}{1.1}
\caption{\textsc{Patient Data From Usecases}}
\label{tab:ucd}
\centering
\begin{tabular}{|l|l|l|l|}
\hline
\multicolumn{1}{|c|}{\textbf{\begin{tabular}[c]{@{}c@{}}Patient Data\end{tabular}}} & \multicolumn{1}{c|}{\textbf{\begin{tabular}[c]{@{}c@{}}Usecase 1\end{tabular}}}                                                       & \multicolumn{1}{c|}{\textbf{\begin{tabular}[c]{@{}c@{}}Usecase 2\end{tabular}}}                                                                                                             & \multicolumn{1}{c|}{\textbf{\begin{tabular}[c]{@{}c@{}}Usecase 3\end{tabular}}}                                                                                      \\ \hline
Respiratory rate                                                                             & 18                                                                                                                                      & 12                                                                                                                                                                                            & 11                                                                                                                                                                   \\ \hline
Blood pressure                                                                               & 145/91                                                                                                                                  & 142/85                                                                                                                                                                                        & 132/82                                                                                                                                                                 \\ \hline
\begin{tabular}[c]{@{}l@{}}Arterial pressure\end{tabular}                           & 104.4                                                                                                                                   & 104                                                                                                                                                                                           & 99                                                                                                                                                                     \\ \hline
Pulse rate                                                                                   & 105                                                                                                                                     & 106                                                                                                                                                                                           & 83                                                                                                                                                                     \\ \hline
Blood glucose                                                                                & 124                                                                                                                                     & 139                                                                                                                                                                                           & 105                                                                                                                                                                     \\ \hline
SpO2                                                                                         & 97                                                                                                                                      & 100                                                                                                                                                                                           & 86                                                                                                                                                                   \\ \hline
Observation                                                                                  & \begin{tabular}[c]{@{}l@{}}Chest pain\\ \\ Respiratory ok\\ \\ No injury\\ \\ Mentally fit\\ \\ Conscious\\ \\ Not alcoholic\end{tabular} & \begin{tabular}[c]{@{}l@{}}No Chest pain\\ \\ Respiratory ok\\ \\ No injury\\ \\ Abdominal pain\\ \\ Pre-abdominal \\ illness\\ \\ Mentally unfit\\ \\ Conscious\\ \\ Not alcohlic\end{tabular} & \begin{tabular}[c]{@{}l@{}}No Chest pain\\  \\ Respiratory low\\ \\ Head injury\\ \\ Mentally unfit\\ \\ Conscious poor\\ \\ Head discomfort\\ \\ Pre-neurological \\ illness\\ \\ Alcoholic\end{tabular} \\ \hline
\end{tabular}
\end{table}
\subsection{Use-Case 1}
Original diagnosis: Unexplained chest pain - Cardiovascular complication. 

The outcome of the XGB and ANN models in Table \ref{tab:uc1} suggests that the patient is experiencing a cardiovascular complication. The patient had chest pain, which was the deciding factor. The second highest probability is for respiratory complications, as respiratory-related chest pain is a common symptom and can suggest a serious or potentially life-threatening condition \cite{brims2010respiratory}.
\begin{table}[!b]
\renewcommand{\arraystretch}{1.2}
\caption{\textsc{Health Complication Probability Prediction For Use Case - 1}}
\label{tab:uc1}
\centering
\begin{tabular}{|l|l|l|l|}
\hline
\multicolumn{1}{|c|}{\textbf{Complication}} & \multicolumn{1}{c|}{\textbf{\begin{tabular}[c]{@{}c@{}}Probability\\ XGB\\ (\%)\end{tabular}}} & \multicolumn{1}{c|}{\textbf{\begin{tabular}[c]{@{}c@{}}Probability\\ Keras\\ (\%)\end{tabular}}} & \multicolumn{1}{c|}{\textbf{\begin{tabular}[c]{@{}c@{}}Probability\\ Pytorch\\ (\%)\end{tabular}}} \\ \hline
Cardiovascular                              & 95.3                                                                                           & 100                                                                                              & 98.94                                                                                              \\ \hline
Respiratory                                 & 70.76                                                                                          & 66.18                                                                                            & 76.81                                                                                              \\ \hline
Neurological                                & 8.15                                                                                           & 13.86                                                                                            & 9.61                                                                                               \\ \hline
Psychiatric                                 & 7.28                                                                                           & 18.8                                                                                             & 18.07                                                                                              \\ \hline
Metabolic                                  & 7.02                                                                                           & 2.92                                                                                             & 12.32                                                                                              \\ \hline
Abdominal                                   & 5.9                                                                                            & 9.93                                                                                             & 15.35                                                                                              \\ \hline
\end{tabular}
\end{table}
\subsection{Use-Case 2}
Original diagnosis: Appendicitis, cyst on the right - Abdominal complication. 

As can be seen in Table \ref{tab:uc2}, the models suggested that the patient had abdominal complications that matched the original diagnosis of the rescue personnel. The deciding factor was the presence of abdominal pain which was diagnosed as appendicitis. The patient also had a history of abdominal disease. The other probabilities were not high enough, except for the outcome of the PyTorch model for metabolic complication. Because, during training period, the model found a close relation of abdominal diseases with metabolic complication. Although, the patient was not fit mentally, the models showed a low probability of psychiatric complication as the other symptoms did not suggest that.
\begin{table}[!t]
\renewcommand{\arraystretch}{1.2}
\caption{\textsc{Health Complication Probability Prediction For Use Case - 2}}
\label{tab:uc2}
\centering
\begin{tabular}{|l|l|l|l|}
\hline
\multicolumn{1}{|c|}{\textbf{Complication}} & \multicolumn{1}{c|}{\textbf{\begin{tabular}[c]{@{}c@{}}Probability\\ XGB\\ (\%)\end{tabular}}} & \multicolumn{1}{c|}{\textbf{\begin{tabular}[c]{@{}c@{}}Probability\\ Keras\\ (\%)\end{tabular}}} & \multicolumn{1}{c|}{\textbf{\begin{tabular}[c]{@{}c@{}}Probability\\ Pytorch\\ (\%)\end{tabular}}} \\ \hline
Abdominal                                   & 83.3                                                                                           & 74.43                                                                                            & 65.94                                                                                              \\ \hline
Cardiovascular                              & 41.09                                                                                          & 28.24                                                                                            & 47.5                                                                                               \\ \hline
Respiratory                                 & 25.41                                                                                          & 27.63                                                                                            & 27.87                                                                                              \\ \hline
Psychiatric                                 & 20.34                                                                                          & 32.23                                                                                            & 34.56                                                                                              \\ \hline
Metabolic                                  & 10.26                                                                                          & 38.41                                                                                            & 63.24                                                                                              \\ \hline
Neurology                                   & 8.15                                                                                           & 13.86                                                                                            & 9.61                                                                                               \\ \hline
\end{tabular}
\end{table}
\subsection{Use-Case 3}
Original diagnosis: Suicidal tendency - Psychiatric complication.
Differential diagnosis: Laceration, Craniocerebral injury - Neurological complication.

The patient in case 3 was mentally ill with poor consciousness. The patient was also found to be alcoholic during rescue showing signs of mental instability. In addition, the patient had an accident that caused a serious head injury. The outcome obtained from Table \ref{tab:uc3} shows that the patient should be treated first for neurological complication, which is justified because the patient complained of head discomfort and the head injury should be prioritized. Moreover, the patient had a history of pre-neurological illness. An important insight into the outcome suggests that although the patient had a low respiratory rate and oxygen saturation, the probability of respiratory complication is not very high. According to Goldman et al. in \cite{goldman2023respiratory}, patients with neurological injuries frequently exhibit irregular breathing patterns that can hinder oxygen from entering the bloodstream.
\begin{table}[!b]
\renewcommand{\arraystretch}{1.2}
\caption{\textsc{Health Complication Probability Prediction For Use Case - 3}}
\label{tab:uc3}
\centering
\begin{tabular}{|l|l|l|l|}
\hline
\multicolumn{1}{|c|}{\textbf{Complication}} & \multicolumn{1}{c|}{\textbf{\begin{tabular}[c]{@{}c@{}}Probability\\ XGB\\ (\%)\end{tabular}}} & \multicolumn{1}{c|}{\textbf{\begin{tabular}[c]{@{}c@{}}Probability\\ Keras\\ (\%)\end{tabular}}} & \multicolumn{1}{c|}{\textbf{\begin{tabular}[c]{@{}c@{}}Probability\\ Pytorch\\ (\%)\end{tabular}}} \\ \hline
Neurology                                   & 99.32                                                                                          & 99.09                                                                                            & 99.6                                                                                               \\ \hline
Psychiatric                                 & 84.29                                                                                          & 96.23                                                                                            & 89.99                                                                                              \\ \hline
Abdominal                                   & 23.77                                                                                          & 22.09                                                                                            & 8.76                                                                                               \\ \hline
Metabollic                                  & 21.78                                                                                          & 52.79                                                                                            & 74.49                                                                                              \\ \hline
Cardiovascular                              & 18.87                                                                                          & 0                                                                                                & 44.38                                                                                              \\ \hline
Respiratory                                 & 11.8                                                                                           & 15.52                                                                                            & 27.87                                                                                              \\ \hline
\end{tabular}
\end{table}
\\
Patient data taken in rescue events often show dynamic characteristics because health vitals are generally not stable and change over time. So for use case-3, the data were modified by a 20\% increase in health vital values. The results can be shown in Table \ref{tab:uc3m} which is similar to the unmodified data indicating that our models were reliable and robust. Based on our research, we also discovered that robustness is attributed to machine learning algorithms that generate output considering not only the individual feature but also the interrelationships between features.
\begin{table}[!t]
\renewcommand{\arraystretch}{1.1}
\caption{\textsc{Health Complication Probability After Modifying Health Vital Data of Use Case - 3}}
\label{tab:uc3m}
\centering
\begin{tabular}{|l|l|l|l|}
\hline
\multicolumn{1}{|c|}{\textbf{Complication}} & \multicolumn{1}{c|}{\textbf{\begin{tabular}[c]{@{}c@{}}Probability\\ XGB\\ (\%)\end{tabular}}} & \multicolumn{1}{c|}{\textbf{\begin{tabular}[c]{@{}c@{}}Probability\\ Keras\\ (\%)\end{tabular}}} & \multicolumn{1}{c|}{\textbf{\begin{tabular}[c]{@{}c@{}}Probability\\ Pytorch\\ (\%)\end{tabular}}} \\ \hline
Neurology                                   & 99.32                                                                                          & 99.09                                                                                            & 99.6                                                                                               \\ \hline
Psychiatric                                 & 80.98                                                                                          & 94.82                                                                                            & 89.78                                                                                              \\ \hline
Abdominal                                   & 5.65                                                                                           & 22.14                                                                                            & 9.99                                                                                               \\ \hline
Metabollic                                  & 11.57                                                                                          & 36.5                                                                                             & 33.45                                                                                              \\ \hline
Cardiovascular                              & 25.46                                                                                          & 0.01                                                                                             & 43.33                                                                                              \\ \hline
Respiratory                                 & 49.02                                                                                          & 14.79                                                                                            & 27.87                                                                                              \\ \hline
\end{tabular}
\end{table}
\section{Conclusion}\label{conc}
The health management software we developed showed a high accuracy for all the use cases and considered multiple aspects in situation recognition. To work with the software in real life events, data collection will be a critical aspect to consider. To function properly, the software requires both health vital data collected from sensors and categorical data to be collected from a set of questionnaire. A comprehensive approach can be to deploy the software in wearable devices which may ease the data collection procedure. Also, the approach used in the software can be further discussed with medical professionals which in turn can improve the decision based outcomes.
\section*{Acknowledgement}
The ongoing research was financially supported by the Federal Ministry of Education and Research, Germany. The research has been supported by the KIRETT project coordinator CRS Medical GmbH (Aßlar, Germany), and partner organization mbeder GmbH (Siegen, Germany).
\bibliographystyle{ieeetr}
\bibliography{reference}

\end{document}